\documentclass[pdflatex,sn-mathphys-num]{sn-jnl}

\usepackage{graphicx}%
\usepackage{multirow}%
\usepackage{amsmath,amssymb,amsfonts}%
\usepackage{amsthm}%
\usepackage{mathrsfs}%
\usepackage[title]{appendix}%
\usepackage{xcolor}%
\usepackage{textcomp}%
\usepackage{manyfoot}%
\usepackage{booktabs}%
\usepackage{algorithm}%
\usepackage{algorithmicx}%
\usepackage{algpseudocode}%
\usepackage{listings}%

\newcommand{\R}[0]{\mathbb{R}}
\newcommand{\N}[0]{\mathbb{N}}

\theoremstyle{thmstyleone}%
\newtheorem{theorem}{Theorem}

\theoremstyle{thmstyletwo}%

\theoremstyle{thmstylethree}%

\begin{document}

\title[]{Cancer Detection via Electrical Impedance Tomography and Optimal Control of Elliptic PDEs}

\author*[1]{\fnm{Ugur G.} \sur{Abdulla}}\email{ugur.abdulla@oist.jp}
\author[1]{\fnm{Jose H.} \sur{Rodrigues}}\email{jose.rodrigues@oist.jp}

\affil*[1]{\orgdiv{Analysis \& PDE Unit}, \orgname{Okinawa Institute of Science and Technology}, \orgaddress{\street{1919-1 Tancha}, \city{Onna-son}, \postcode{904-0495}, \state{Okinawa}, \country{Japan}}}

\abstract{
We pursue a computational analysis of the biomedical problem on the identification of the cancerous tumor at an early stage of development based on the Electrical Impedance Tomography (EIT) and optimal control of elliptic partial differential equations.
Relying on the fact that the electrical conductivity of the cancerous tumor is significantly higher than the conductivity of the healthy tissue, we consider
an inverse EIT problem on the identification of the conductivity map in the complete electrode model based on the $m$ current-to-voltage measurements on the boundary electrodes. 
A variational formulation as a PDE-constrained optimal control problem is introduced based on the novel idea of increasing the size of the input data by adding "voltage-to-current" measurements through various permutations of the single "current-to-voltage" measurement. The idea of permutation preserves the size of the unknown parameters on the expense of increase of the number of PDE constraints. We apply a gradient projection (GPM) method based on the Fr\'echet differentiability in Besov-Hilbert spaces. 
Numerical simulations of 2D and 3D model examples demonstrate the sharp increase of the resolution of the cancerous tumor by increasing the number of measurements from $m$ to $m^2$.\\
}



\maketitle

\section{Background}\label{intro}

This paper addresses the inverse EIT problem of detecting an unknown conductivity inside a body, based on voltage measurements on the surface when electric currents are applied through a finite set of electrodes.
Let $Q\subset\R^n$ be an open and bounded domain, and $\sigma:Q\to\mathbb{R}$ be the conductivity map. Let $E=(E_l)_{l=1}^m$, $m\in\N$, be a finite set of electrodes attached to the surface $\partial Q$, with corresponding \textit{contact impedance} vector $Z=(Z_l)_{l=1}^m\in\R^m$. The electric \textit{current pattern} vector $I=(I_l)_{l=1}^m\in\R^m$ is applied to the electrodes $E$ and the corresponding induced \textit{voltage} vector $U=(U_l)_{l=1}^m\in\R^m$ is measured. The following {\it conservation of charges} and \textit{grounding conditions} are satisfied:
\begin{equation*}
    \sum_{l=1}^m I_l = 0 \ \  \quad \sum_{l=1}^m U_l = 0.
\end{equation*}
The potential $u$ inside the body $Q$ is described by the following second order elliptic partial differential equation and corresponding boundary conditions:
\begin{eqnarray}
    -\mbox{div}(\sigma(x)\nabla{u}) = 0, &\mbox{in}& Q \label{pde}\\
    \sigma(x)\frac{\partial u}{\partial\nu} = 0, &\mbox{on}& \partial Q - \bigcup_{l=1}^m E_l \label{bc_1}\\
    u + Z_l \sigma(x)\frac{\partial u}{\partial\nu} = U_l, &\mbox{on}& E_l,~l=1,\dots,m \label{bc_2}\\
    \int_{E_l}\frac{\partial u}{\partial\nu}ds = I_l, & & l=1,\dots,m \label{bc_3}
\end{eqnarray}
where $\nu(x)=(\nu^i(x))_{i=1,\dots,n}$ is the outward normal at the point $x\in\partial Q$.\\

{\bf Inverse EIT Problem:} {\it Given electrode contact impedance vector $Z$, electrode {\it current pattern} $I$ and boundary electrode measurement $U^*$, it is required to find electrostatic potential $u$ and electrical conductivity map $\sigma$ satisfying \eqref{pde}--\eqref{bc_3}} with $U=U^*$.

The inverse EIT problem is highly ill-posed, as it aims to identify an infinite-dimensional conductivity map using finite-dimensional "current-to-voltage" measurements on the electrodes. Recently, a new variational method has been introduced in \cite{Abdulla2021} based on the PDE constrained optimal control problem in Sobolev space setting. The novelty of the control-theoretic model is its adaptation to the clinical situation when additional "voltage-to-current" measurements based on the various permutations of the single "current-to-voltage" measurement can increase the size of the input data while keeping the size of the unknown parameters fixed. In \cite{Abdulla2021} existence of the optimal control and Fr\'echet differentiability in Besov-Hilbert spaces is proved, the formula for the Fr\'echet gradient is derived and a gradient descent algorithm in Besov-Hilbert spaces has been developed. In \cite{Abdulla2023} convergence of the finite-difference method is proved. The goal of this paper is to develop a computational framework based on the gradient projection method in the Besov-Hilbert spaces to identify cancerous tumors both in 2D and 3D model examples. 

EIT problems have a large number of applications in medicine, industry, geophysics and material sciences \cite{Holder2004}. Forward EIT problem \eqref{pde}-\eqref{bc_3} for the identification of $(u,U)$ with given input data $(\sigma, I, Z)$  is referred as \textit{complete electrode model}. It was introduced in \cite{Somersalo1992} as a physically more accurate model capable of predicting experimental data with high precision. Existence and uniqueness of solution to the complete electrode model \eqref{pde}-\eqref{bc_3} was established in \cite{Somersalo1992}.

Motivated by the medical applications on detection of cancerous tumors from the breast tissue and other parts of the body, the relevance of inverse EIT problem resides on the fact that the conductivity of cancerous tissue is considerably higher than the conductivity of normal tissue \cite{Halter2009,Laufer2010}.

The inverse EIT problem belongs to the class of so called Calderon type ill-posed inverse problems due to the celebrated work \cite{Calderon1980}, where the well-posedness of the inverse problem for the identification of $\sigma$ through Dirichlet-to-Neumann or Neumann-to-Dirichlet maps for the PDE \eqref{pde} is presented. Significant development in Calderon's inverse problem concerning questions on uniqueness and stability was achieved in \cite{Sylvester1987,Nachman1988,Kenig2007,Nachman1995,Knudsen2008}. 

The difficulty in solving the inverse EIT problem is due to the identification of the infinite-dimensional conductivity map $\sigma$ and the finite-dimensional voltage vector $U$ using the finitely many measurements of input data. It is important to notice that the number of input data depends on the number of electrodes and there is no flexibility to increase its size. Alternatively, an increase of measurement sets (current patterns) could be used to identify the same conductivity map, however, the number of unknown voltages would increase accordingly. A variety of numerical methods have been developed in the attempt to solve the inverse EIT problem \cite{Ammari2004,Kwon2002,Widlak2012,Seo2011,Alberti2016,Gehre2014,Jin2017,Siltanen2001,Alsaker2016,Dodd2014,Hamilton2016,Hyvonen2018}. 

The majority of the methods mentioned above and reconstruction algorithms found in the literature are dedicated to 2D inverse EIT problem. Therefore, it is natural to expect that 2D algorithms and methods could be used in an attempt to identify anomalies in a cross-section of a 3D body. However, 3D characteristics of the current flow may be neglected creating distortions in the resulting images \cite{Blue2000}. Similar drawback has been previously observed in \cite{Kleinermann1996} and a reconstruction algorithm based on the inversion of the sensitivity matrix was proposed for a simplified model in a finite right circular cylinder. 

\section{Methods} \label{sct:MM}

\subsection{Optimal Control Problem}\label{sct:OCP}
We aim to formulate an inverse EIT problem as an optimal control problem by selecting conductivity map $\sigma$ and boundary electrode voltage vector $U$ as control parameters.
Given control vector $(\sigma,U)$, the state vector-potential $u$ is identified as a Sobolev-Hilbert solution of the elliptic PDE problem \eqref{pde}-\eqref{bc_2}. Optimal control framework is implemented to identify the pair $(\sigma,U)$ which is the best candidate to fulfill the Ohm's law on the electrodes (condition \eqref{bc_3}). Particular advantage of this approach is a well-posedness of the elliptic problem \eqref{pde}-\eqref{bc_2} under very general assumptions on the conductivity map $\sigma$ as a consequence of the powerful Lax-Milgram theory. 

For a given $v=(\sigma,U)\in L_\infty(Q)\times\R^m$, a function $u=u(~\cdot~;v)\in H^1(Q)$ is called a \textit{solution} to the PDE problem \eqref{pde}-\eqref{bc_2} if the following identity is satisfied:
\begin{equation}
    B[u,\eta] = \sum_{l=1}^m\frac{U_l}{Z_l}\int_{E_l}\eta ds,\quad \forall\eta\in H^1(Q).
\end{equation}

To prove the necessary optimality condition we introduce the adjoined state problem corresponding to \eqref{pde}-\eqref{bc_2}. Given a control vector
$v=(\sigma,U)\in L_\infty(Q)\times\R^m$, let $u=u(~\cdot~;v)$ be the corresponding solution of \eqref{pde}-\eqref{bc_2}. The following is the \textit{adjoined} problem to \eqref{pde}-\eqref{bc_2}. 
\begin{eqnarray}
    -\mbox{div}(\sigma(x)\nabla\psi) = 0, &\mbox{in}& Q \label{adjn_pde}\\
    \sigma(x)\frac{\partial\psi}{\partial\nu} = 0, &\mbox{on}& \partial Q - \bigcup_{l=1}^m E_l \label{adjn_bc1}\\
    \psi + Z_l\sigma(x)\frac{\partial\psi}{\partial\nu} = 2\int_{E_l}\frac{u(s)-U_l}{Z_l}dS + 2I_l, &\mbox{on}& E_l,\quad l=1,\dots,m . \label{adjn_bc2}
\end{eqnarray}
A function $\psi=\psi(~\cdot~;v)\in H^1(Q)$ is a \textit{solution} to \eqref{adjn_pde}-\eqref{adjn_bc2} if the following identity is verified
\begin{equation}
    B[\psi,\eta] = \sum_{l=1}^m\frac{2}{Z_l}\left(\int_{E_l}\frac{u-U_l}{Z_l}ds+I_l\right)\int_{E_l}\eta ds,\quad\forall\eta\in~H^1(Q).
\end{equation}
The existence, uniqueness and stability results for the solutions to elliptic PDE problems \eqref{pde}-\eqref{bc_2} and \eqref{adjn_pde}-\eqref{adjn_bc2} is a consequence of Lax-Milgram theory in Sobolev-Hilbert space $H^1(Q)$ \cite{Abdulla2021}. 
 
Consider the following variational formulation of the inverse EIT Problem: given electrode current pattern $I$ and corresponding electrode voltage measurement vector $U^*$, consider the minimization of the functional
\begin{equation}
    \mathcal{J}(v) = \sum_{l=1}^{m}\left| \int_{E_l}\frac{U_l-u}{Z_l} -I_l\right|^2 + \beta|U-U^\ast|^2, \quad \beta>0, \label{cost_j}
\end{equation}
on the control set
\begin{gather}
V_R=\left\{
\begin{array}{c}
v = (\sigma, U)\in \Big( L_{\infty} (Q) \cap H^{\epsilon}(Q) \Big )
\times \mathbb{R}^m \Big| \\ 
\displaystyle\sum_{l=1}^m U_l = 0,
\|\sigma\|_{L_{\infty}}+ \|\sigma\|_{H^{\epsilon}} +|U| \leq R, \  \sigma \geq \mu > 0  
\end{array}
\right\}\nonumber
\end{gather}
where $\beta> 0$, and $u=u(~\cdot~;v)\in H^1(Q)$ is the solution of \eqref{pde}-\eqref{bc_2}.
This optimal control problem will be called Problem $\mathcal{J}$. The first term in the cost functional \eqref{cost_j} represents the error for integral from of the Ohm's law on the boundary electrodes (condition \eqref{bc_3}) in light of the Robin boundary condition \eqref{bc_2}.

It should be stressed out that the variational formulation of the forward EIT Problem is a particular case of the Problem $\mathcal{J}$. If the conductivity map $\sigma$ is known, we consider the optimal control problem on the minimization of the function 
\begin{align}
  \mathcal{I}(U) =   \displaystyle\sum _{l=1}^m \Big |  \displaystyle\int_{E_l} \frac{U_l-u(x)}{Z_l}ds- I_l\Big|^2\label{eq:cost_functional1}\to \inf
\end{align}
in a control set
\begin{equation}\label{convex2}
W=\big\{U\in \mathbb{R}^m \Big| \sum_{l=1}^m U_l = 0 \big\}
\end{equation}
where $u=u(\cdot;v)\in H^1(Q)$ is a solution of the elliptic problem \eqref{pde}--\eqref{bc_2}.
Furthermore this optimal control problem will be called Problem $\mathcal{I}$. It is a convex PDE constrained optimal control problem.

The optimal control problem $\mathcal{J}$ inherits the ill-posedness of the inverse EIT problem. Following \cite{Abdulla2021} we now formulate the optimal control problem which is adapted to the situation when the size of the input data can be increased 
through additional measurements while keeping the size of the unknown parameters fixed. Let $I^1:=I$ is a current pattern input, and $U^1=(U_1,...,U_m)$ is a corresponding boundary electrode voltage measurement. Consider
$m-1$  new permutations of boundary voltages
\begin{equation}
U^j = (U_j,\dots,U_m,U_1,\dots,U_{j-1}), \quad j=2,\dots,m, \label{rot_schm}
\end{equation}
of $U$. The set of permutations above will be referred to as "rotation scheme". Assume that the "voltage-to-current" measurement allows us to measure associated current $I^j=(I^j_1,...,I^j_m)$, $j=1,\dots,m$. By setting $U^1=U^*$ and having a new set of $m^2$ input data $(I^j)^{m}_{j=1}$,
we now consider optimal control problem on the minimization of the new cost functional
\begin{equation}
    \mathcal{K}(v) = \sum_{j=1}^{m}\sum_{l=1}^m\left|\int_{E_l}\frac{U^j_l-u^j}{Z_l}dS - I^j_l \right|^2 + \beta|U-U^\ast|^2,
\end{equation}
on a control set $V_R$, where $u^j=u(~\cdot~;\sigma,U^j)$ is the corresponding solution to the problem \eqref{pde}-\eqref{bc_2} with $U$ replaced by its permutation $U^j$, for $j=1,\dots,m$. This optimal control problem will be called Problem $\mathcal{K}$. Note that the number of input currents in the Problem $\mathcal{K}$ has increased from $m$ to 
$m^2$. However, the size of unknown control vector is unchanged, and in particular there are only $m$ unknown voltages $U_1,\cdots,U_m$, whereas all vectors $U^j, j=2,...,m$ are formed by their permutation as in \eqref{rot_schm}. The price we pay for this gain is the increase of the number of PDE constraints from 1 to $m$. It is essential to note that the 
Problem $\mathcal{J}$ is a particular case of the Problem $\mathcal{K}$, precisely when we don't consider any permutation of $U^1$, but the trivial one.

\subsection{Fr\'echet Differentiability}
Existence of an optimal control for the problem $\mathcal{K}$ (consequently $\mathcal{J}$) and Fr\'echet differentiability was proved in \cite{Abdulla2021}.
In \cite{Abdulla2023} the convergence of the method of finite differences is established. 

\begin{theorem}[Fr\'echet Differentiability] (\cite{Abdulla2021}): 
    The functional $\mathcal{K}$ is differentiable on $V_R$ in the sense of Fr\'echet and the corresponding Fr\'echet gradient $\mathcal{K}':V_R\to\textbf{ba} (Q) \times\R^m$ is given by $\mathcal{K}'(v) = \left(\mathcal{K}'_\sigma(\sigma,U),\mathcal{K}'_U(\sigma,U)\right)$ where
    \begin{eqnarray}
        \mathcal{K}'_\sigma(\sigma,U) &=& -\sum_{j=1}^{m}\nabla \psi^j \cdot \nabla u^j, \label{grad_A}\\
        \mathcal{K}'_U(\sigma,U) &=& \left(
               \sum_{j,l=1}^{m} 2\left[\int_{E_l}\frac{U^j_l-u^j(s)}{Z_l}dS-I^j_l\right]\int_{E_l}\frac{\delta_{l\theta_{kj}}-w^{\theta_{kj}}(s)}{Z_l}dS
                +2\beta(U_k-U^\ast_k) 
            \right)_{k=1}^m . \label{grad_U}
    \end{eqnarray}
    where
    $$ \theta_{kj} = \left\{\begin{array}{ll}
        k-j+1, & \mbox{if}~ j \leq k, \\
        m+k-j+1, & \mbox{if}~ j > k,
    \end{array}\right. $$
    and $w^{\theta_{kj}}=u(~\cdot~;\sigma,e_{\theta_{kj}})$ and $e_{\theta_{kj}}\in\R^m$ is the unit vector in the $\theta_{kj}$-direction.
\end{theorem}

\subsection{Gradient Method in Besov-Hilbert Space} \label{ssct:grad_mtd}
Fr\'echet differentiability result suggest the following algorithm based on the projective gradient method for the Problem $\mathcal{K}$. \\
\noindent\textbf{Step 1.} 
Set iteration counter $N=0$ and choose initial controls $v^0=(\sigma^0,U^0)\in V_R$, where $U^0=(U^0_l)_{l=1}^m$ satisfies $\sum_{l=1}^m U^0_l = 0$.

\vskip.1in

\noindent\textbf{Step 2.} 
Build the pertmutations $U^{N,j}$, solve the problem \eqref{pde}-\eqref{bc_2} to find $u^{N,j}=u(~\cdot~;\sigma^N,U^{N,j})$, $j=1,\dots,m$, and calculate $\mathcal{K}^N = \mathcal{K}(\sigma^N,U^N)$.

\vskip.1in

\noindent\textbf{Step 3.} 
If $N=0$, move to Step 4. Otherwise, check for the error condition
\begin{equation}
    \max\left( 
    \left|\frac{\mathcal{K}^N-\mathcal{K}^{N-1}}{\mathcal{K}^{N-1}}\right|,
    \frac{|U^N-U^{N-1}|}{|U^{N-1}|},
    \frac{\|\sigma^N-\sigma^{N-1}\|_{L_2}}{\|\sigma^{N-1}\|_{L_2}}
    \right) < \epsilon, \label{rel_err}
\end{equation}
where $\epsilon>0$ is the required accuracy. If \eqref{rel_err} is verified, then terminate the iteration process. Otherwise, move to Step 4.

\vskip.1in

\noindent\textbf{Step 4.} 
Solve the problem \eqref{pde}-\eqref{bc_2} to find $w^N_k=u(~\cdot~;\sigma^N;e_k)$, where $e_k\in\R^m$ is the unit vector in the $k$-direction, $k=1,\dots,m$.

\vskip.1in

\noindent\textbf{Step 5.} 
Solve the adjoined problem \eqref{adjn_pde}-\eqref{adjn_bc2} to find adjoined potential $\psi^{N,j}=\psi(~\cdot~;\sigma^N,U^{N,j},u^{N,j})$, for $j=1,\dots,m$.

\vskip.1in

\noindent\textbf{Step 6.}
Choose step size parameter $\gamma^N>0$ and compute new control vector components $v^{N+1}=(\sigma^{N+1},U^{N+1})$ using the Fr\'echet differentiability expressions \eqref{grad_A} and \eqref{grad_U} as follows
\begin{align}
    \tilde{\sigma}^{N+1}(x) &= \sigma^N(x) + \gamma^N
    \sum_{j=1}^{m}\nabla \psi^{N,j}(x) \cdot \nabla u^{N,j}(x),\quad x\in Q, \\
    \tilde{U}^{N+1}_k &= U^{N+1}_k \\
    	& - \gamma^N \left[ 
        \sum_{j,l=1}^{m} 2\left( \int_{E_l}\frac{U^{N,j}_{l}-u^{N,j}}{Z_l}dS \right) \int_{E_l}\frac{\delta_{l,\theta_{kj}}-w^{\theta_{kj}}}{Z_l}dS + 2\beta(U^N_k-U^\ast_k)
        \right] ,
\end{align}
for $k=1,\dots,m$.

\vskip.1in

\noindent\textbf{Step 7.}
Update new control components using the following projection formulas
\begin{eqnarray}
    \sigma^{N+1}(x) &=& \left\{\begin{array}{ll}
        \mu, & \mbox{if}\quad\tilde{\sigma}^{N+1}(x) \leq \mu, \\
        \tilde{\sigma}^{N+1}(x), & \mbox{if}\quad\mu \leq \tilde{\sigma}^{N+1}(x) \leq R, \\
        R, & \mbox{if}\quad\tilde{\sigma}^{N+1}(x) \geq R,
    \end{array}\right. , \quad x\in Q, \\
    U^{N+1}_k &=& \tilde{U}^{N+1}_k - \frac{1}{m}\sum_{l=1}^{m}\tilde{U}^{N+1}_l,\quad k=1,\dots,m.
\end{eqnarray}
Then, replace $N$ by $N+1$ and move to Step 2. \\

\subsection{Two-stage Procedure for Simulations and Clinical Applications.}\label{multistage}
The above algorithm suggest the following two-stage procedure with the increase of data both in simulations, as well as clinical applications.\\

\begin{itemize}
\item {\it Simulation-Stage 1}: Select one set of electrode current input $I^1=(I^1_l)_{l=1}^m$ and the conductivity map $\sigma=\sigma_{true}$ reflecting a distribution of cancerous tumor; solve the convex optimization Problem $\mathcal{I}$ and find its unique minimizer $U_{true}$. Let $U^*=U_{true}$. The pair $\left( \sigma_{true}, u_{true} \right)$ is the solution of the inverse EIT problem with input data $(Z,I^1,U^*)$. Here, $u_{true}=u(\cdot; \sigma_{true}, U^*)$ be a solution of the elliptic PDE problem \eqref{pde}-\eqref{bc_2}. 
\item {\it Simulation-Stage 2}: Denote $U^1=U^*$ and consider $m-1$ permutations $\{U^j\}_{j=2}^m$ as in \eqref{rot_schm}. Denote $u^1\equiv u_{true}$, and for each $j=2,...,m$, solve the elliptic PDE problem  \eqref{pde}-\eqref{bc_2} to find functions $u^j=u(\cdot; \sigma_{true}, U^j)$. Then use the "voltage-to-current" formula \eqref{bc_3} with $u=u^j$ to calculate $m-1$ new sets of current vectors $I^j, j=2,...,m$. This procedure guarantees that $\left( \sigma_{true}, U_{true} \right)$ is an optimal control for the Problem $\mathcal{K}$. Solve the optimal control Problem $\mathcal{K}$ with $m^2$ input data $I=(I^j_l)_{j,l=1}^m$ by the GPM algorithm to recover an optimal control $\left( \sigma_{true}, U_{true} \right)$.
\end{itemize}

Our main results suggest the following two-stage procedure for medical application for the identification of the cancerous tumor at an early stage of development:
\begin{itemize}
\item \textit{Clinical Application-Step 1}: Apply selected one set of electrode current vector $I^1=(I^1_l)_{l=1}^m$ on the electrodes $E=(E_l)_{l=1}^m$, take the "current-to-voltage" measurements $U^1=(U^1_l)_{l=1}^m$.
\item \textit{Clinical Application-Step 2}: Consider $m-1$ permutations $\{U^j\}_{j=2}^m$ as in \eqref{rot_schm}; apply each voltage vectoe $U^j$ to electrodes; pursue "voltage-to-current" measurements $I^j=(I^j_l)_{l=1}^m$, and then solve the optimal control Problem $\mathcal{K}$ with $m^2$ input data $I=(I^j_l)_{j,l=1}^m$ by the GPM algorithm to identify the location of development of the cancerous tumor. 
\end{itemize}

\section{Results}\label{sct:Res}

\subsection*{Methodology}

\noindent\textit{Finite element approach.}
The PDE problem \eqref{pde}-\eqref{bc_2} and its adjoined \eqref{adjn_pde}-\eqref{adjn_bc2} are numerically solved using the \textit{partial differential equation toolbox} package of Matlab. This package applies a linear \textit{finite element} spatial discretization of the domain $Q=\bigcup Q^e$, where $Q^e$ is a tetrahedron element. If $\{\phi_i\}_{i=1,\dots,N_d}$ denotes the piecewise polynomial basis of functions, where $N_d$ is the number of nodes in the discretization, then the approximate solution of problem \eqref{pde}-\eqref{bc_2} is written as
$u_{\mbox{fem}} = \sum_{i=1}^{N_d}u^i\phi_i$
where $u^i$ is the undetermined scalar corresponding to the potential $u$ at the node $i$. Hence, the variational formulation yields the linear system
\begin{equation}
    \sum_{i=1}^{N_d}u^i
    \left[\int_Q \sigma \nabla{\phi_i}\cdot\nabla{\phi_j}dx
        + \sum_{l=1}^{m}\frac{1}{Z_l}\int_{E_l}\phi_i\phi_jdS \right]
    = \sum_{l=1}^{m}\frac{U_l}{Z_l}\int_{E_l}\phi_jdS,\quad j=1,\dots,N_d. \label{num_frwd_var}
\end{equation}
Analogously, if 
$\psi_{\mbox{fem}}=\sum_{i=1}^{N_d}\psi^i\phi_i$ corresponds to the discretization of the adjoined potential $\psi$, then the corresponding variational formulation yields
\begin{align}
    \sum_{i=1}^{N_d}\psi^i
    \left[\int_Q \sigma \nabla{\phi_i}\cdot\nabla{\phi_j}dx
        + \sum_{l=1}^{m}\frac{1}{Z_l}\int_{E_l}\phi_i\phi_jdS \right] \nonumber \\
    = \sum_{l=1}^{m}\int_{E_l}\frac{\phi_j}{Z_l}\left[2\int_{E_l}\frac{u_{\mbox{fem}}-U_l}{Z_l}dS+2I_l\right]dS, \label{num_bwrd_var}
\end{align}
for $j=1,\dots,N_d$.

For the inverse EIT problem, we start by setting current $I$ and contact impedance $Z$ vectors as described in the previous section. In order to simulate the EIT model, we set the conductivity map $\sigma_{true}:Q\to\R$ to emulate spherical tumor regions of center $c$ and radius $r>0$ inside $Q$, namely
\begin{equation}
    \sigma_{true}(x) = 
    \left\{\begin{array}{ll}
        0.4, &\mbox{if}~|x-c|\leq r; \\
        0.2, &\mbox{otherwise},
    \end{array}\right. \qquad\mbox{in}~(\mbox{Ohm}\cdot\mbox{m})^{-1} .\label{sigma_true}
\end{equation}

\noindent\textit{Choice of learning rate parameter $\alpha$.}
For all the simulations, the learning rate parameter $\alpha^N$ in Step 6 of the GPM algorithm described in Section 4.1 of \cite{Abdulla2021} was calculated in each iteration as the average of Barzilai-Borwein -type formulas \cite{Barzilai1988}. Indeed, separate coefficients were calculated for each variable (voltage and conductivity) as follows:
\begin{eqnarray}
    \alpha^N_U = \mbox{mean}\left( 
    \frac{|dU^N|^2}{|dU^N \cdot d\mathcal{K'}_U^N|} ~,~
    \frac{|dU^N \cdot d\mathcal{K'}_U^N|}{|d\mathcal{K'}_U^N|^2}
    \right),\quad \nonumber\\
    \alpha^N_\sigma = \mbox{mean}\left(
    \frac{\|d\sigma^N\|_{L_2}^2}{(d\sigma^N,d\mathcal{K'}_\sigma^N)_{L_2}} ~,~
    \frac{(d\sigma^N,d\mathcal{K'}_\sigma^N)_{L_2}}{\|d\mathcal{K'}_\sigma^N\|_{L_2}^2}
    \right)
\end{eqnarray}
where
$$dU^N = U^{N}-U^{N-1},\quad d\mathcal{K'}_U^N = \mathcal{K'}_U(\sigma^N,U^N)-\mathcal{K'}_U(\sigma^{N-1},U^{N-1}),$$
$$d\sigma^N = \sigma^{N}-\sigma^{N-1},\quad d\mathcal{K'}_\sigma^N = \mathcal{K'}_\sigma(\sigma^N,U^N)-\mathcal{K'}_\sigma(\sigma^{N-1},U^{N-1}).$$

\subsection{Results in 2D}  
\begin{figure}[h!]
\centering
\includegraphics[width=0.9\textwidth]{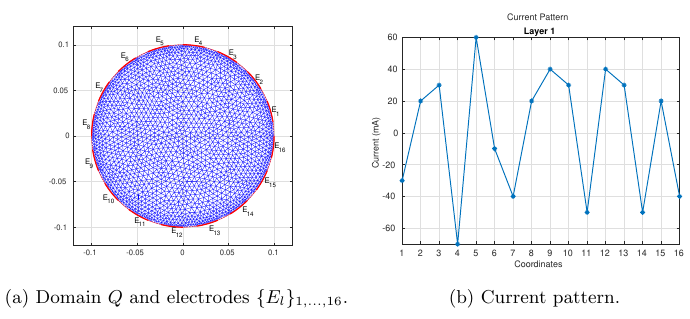}
\caption{(Left) Meshed domain $Q$ and electrodes position. (Right) Current pattern used for all cases.}
\label{fig_2D_01}
\end{figure}

In this case, we set $Q$ to be the circle of radius $r=0.1\mbox{m}$ given by
\begin{equation*}
    Q = \left\{ (x,y) \in \mathbb{R}^2 :~ x^2+y^2<r^2 \right\}.
\end{equation*}
A set of $m=16$ electrodes with dimension $0.024~\mbox{rad}$ width were uniformly distributed along the boundary $\partial Q$, see Figure \ref{fig_2D_01}(a). A mesh consisting of 2034 nodes and 3794 linear elements (triangles) was considered. A uniform contact impedance vector $Z=(Z_l)_{l=1,\dots,16}$ with $Z_l=0.1~\mbox{Ohm}$ was set. Background conductivity is set to $0.2~(\mbox{Ohm}\cdot\mbox{m})^{-1}$ corresponding to healthy tissue, while tumorous tissue corresponds to $0.4~(\mbox{Ohm}\cdot\mbox{m})^{-1}$. The current vector $I$ is set as shown in Figure \ref{fig_2D_01}(b).

For the optimization process using the GPM algorithm described in Section 4.1 in \cite{Abdulla2021}, we set initial conductivity map $\sigma_{ini} = 0.3~(\mbox{Ohm}\cdot\mbox{m})^{-1}$ and initial voltage vector $U^{ini} = (U^{ini}_{l})_{l=1,\dots,m}$ as follows: $U^{ini}_l=1~\mbox{volt}$ if $l$ is even, $U^{ini}_{l}=-1~\mbox{volt}$ if $l$ is odd. Finally, termination conditions were set to a max number of iteration $N_{max}=250$ or relative error tolerance $\epsilon=10^{-6}$. 

\subsubsection{Case: 1 Tumor}
\begin{figure}[h!]
\centering
\includegraphics[width=0.9\textwidth]{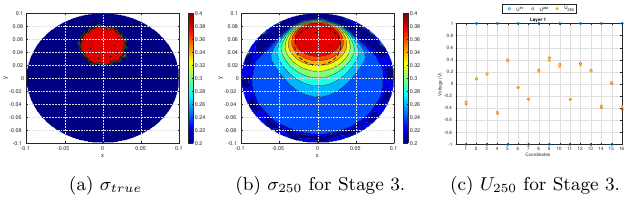}
\caption{Contour plots. (Left) True conductivity $\sigma_{true}$. (Center) Obtained conductivity $\sigma_{250}$ for stage 3. (Right) Obtained voltage coordinates for stage 3.}
\label{fig_2D1T_01}
\end{figure}

We first consider the case of $\sigma_{true}$ with center $c=(0,-0.05)$ and radius $r=0.03$ see Figure \ref{fig_2D1T_01}(a). Optimal control framework is implemented without regularization $(\beta=0)$. Figure \ref{fig_2D1T_01}(b) shows $\sigma_N$ at final iteration $N=250$ for stage 2. Dashed lines show the position and size of the target $\sigma_{true}$. Figure \ref{fig_2D1T_01}(c) shows the coordinates of optimal voltage $U^\ast$, initial voltage $U^{ini}$ and obtained voltage $U^{N}$ at the last iteration $N=250$ of stage 3. Cost value at the final iteration of stage 3 is $\mathcal{K}_{end} = 3.1588\mbox{e-07}$ and relative errors of voltage and conductivity are $\frac{|U^{end}-U^\ast|}{|U_\ast|} = 0.0787$ and $\frac{\|\sigma_{end}-\sigma_{true}\|_{L_2}}{\|\sigma_{true}\|_{L_2}} = 0.2757$.

\begin{figure}[h!]
\centering
\includegraphics[width=0.9\textwidth]{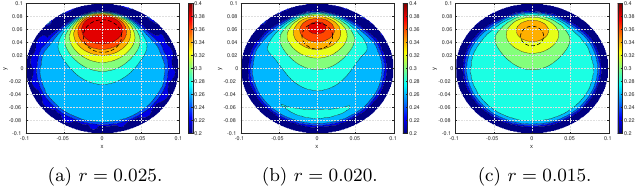}
\caption{Contour plot of obtained conductivity $\sigma_{250}$ for different values of radius $r$ of true conductivity. Dashed lines represent size and position of true conductivity.}
\label{fig_2D1T_02}
\end{figure}

\vskip.1in
\noindent\textbf{Sensitivity with respect to size.}
In order to demonstrate the sensitivity of the method and calculations, we fix the position $c$ of $\sigma_{true}$ and consider different values of radius, namely $r=0.025,0.020,0.015,\allowbreak0.010,0.005$. Figure \ref{fig_2D1T_02} shows the reconstructed conductivity $\sigma_{N}$, at iteration $N=250$, for case $r=0.25,0.20,0.15$. Dashed lines show the position and size of $\sigma_{true}$. Table \ref{tab1} shows the cost values and relative error of voltage and conductivity at the last iteration of stage 2 for each case of radius $r$.

\begin{table}[h]
\caption{Metrics for the 2D - Case: 1 tumor - Size.\label{tab1}}
\begin{tabular}{cccc}
\toprule
\multirow{2}{5em}{Radius ($r$)} & \multirow{2}{5em}{Cost Value ($\mathcal{K}_{end}$)} & \multicolumn{2}{c}{Relative Error}  \\
& & Voltage ($\frac{|U^{end}-U^\ast|}{|U^\ast|}$) & Conductivity ($\frac{\|\sigma_{end}-\sigma_{true}\|_{L_2}}{\|\sigma_{true}\|_{L_2}}$) \\
\midrule
0.030 & 5.2984\mbox{e-07} & 0.0787 & 0.2757 \\
0.025 & 3.2890\mbox{e-07} & 0.0830 & 0.3406 \\
0.020 & 1.6377\mbox{e-07} & 0.0874 & 0.3642 \\
0.015 & 6.6754\mbox{e-08} & 0.0917 & 0.3907 \\
0.010 & 8.8539\mbox{e-09} & 0.0946 & 0.4051 \\
0.005 & 7.2273\mbox{e-10} & 0.0960 & 0.4110 \\
\bottomrule
\end{tabular}
\end{table}

\begin{figure}[h!]
\centering
\includegraphics[width=0.9\textwidth]{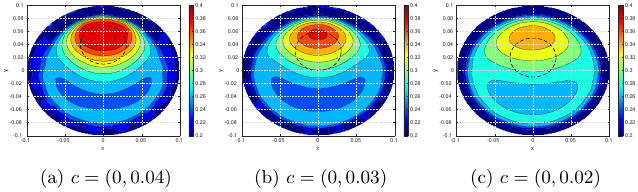}
\caption{Contour plot of obtained conductivity $\sigma_{250}$ for different values of center $c$ of true conductivity. Dashed lines represent size and position of true conductivity.}
\label{fig_2D1T_03}
\end{figure}

\vskip.1in
\noindent\textbf{Sensitivity with respect to position.}
For this analysis, we fixed the radius $r$ of $\sigma_{true}$ and considered different positions of center $c=(0,y)$, $y=0.05,0.04,0.03,0.02,0.01$. Figure \ref{fig_2D1T_03} shows the reconstructed conductivity $\sigma_{N}$ at iteration $N=250$ for cases $y=0.04,0.03,0.02$. Dashed lines show the location of $\sigma_{true}$. Table \ref{tab2} shows the cost value and relative error of voltage and conductivity at the last iteration for each case of center $c$.

\begin{table}[h]
\caption{Metrics for the 2D - Case: 1 tumor - Location.\label{tab2}}
\begin{tabular}{cccc}
\toprule
\multirow{2}{5em}{Center ($c=(0,y)$)} & \multirow{2}{5em}{Cost Value ($\mathcal{K}_{end}$)} & \multicolumn{2}{c}{Relative Error}  \\
& & Voltage ($\frac{|U^{end}-U^\ast|}{|U^\ast|}$) & Conductivity ($\frac{\|\sigma_{end}-\sigma_{true}\|_{L_2}}{\|\sigma_{true}\|_{L_2}}$) \\
\midrule
$0.05$ & 5.2984\mbox{e-07} & 0.0787 & 0.2757 \\
$0.04$ & 3.9891\mbox{e-07} & 0.0776 & 0.3119 \\
$0.03$ & 1.7854\mbox{e-07} & 0.0797 & 0.3089 \\
$0.02$ & 1.1384\mbox{e-07} & 0.0791 & 0.3480 \\
$0.01$ & 2.7750\mbox{e-08} & 0.0795 & 0.3582 \\
$0.00$ & 6.0320\mbox{e-11} & 0.0794 & 0.3615 \\
\bottomrule
\end{tabular}
\end{table}

\subsubsection{Case: 4 Tumors}
\begin{figure}[h!]
\centering
\includegraphics[width=0.9\textwidth]{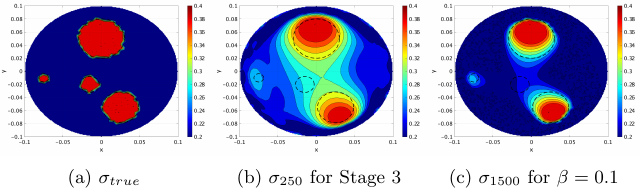}
\caption{(Left) True conductivity showing position and size of 4 tumors. (Center) Contour plot of obtained conductivity $\sigma_{250}$ for stage 3. (Right) Contour plot of obtained conductivity with regularization parameter $\beta=0.1$ after 1500 iterations of the Gradient Method.}
\label{fig_2D1T_04}
\end{figure}

We consider here $\sigma_{true}$ describing 4 circular tumor regions given by centers $c_1=(0,0.050)$, $c_2=(0.025,-0.055)$, $c_3=(-0.015,-0.020)$, $c_4=(-0.075,-0.010)$, and corresponding radius $r_1=0.03$, $r_2=0.0235$, $r_3=0.0122$ and $r_4=0.0063$. True conductivity $\sigma_{true}$ is shown in Figure \ref{fig_2D1T_04}(a). Figure \ref{fig_2D1T_04}(b) shows the reconstructed conductivity for stage 3. The cost value at the final iteration of stage 2 is $\mathcal{K}_{end}=5.5801\mbox{e-07}$ and the corresponding relative errors of voltage and conductivity are $\frac{|U^{end}-U^\ast|}{|U^\ast|}=0.0610$ and $\frac{\|\sigma_{end}-\sigma_{true}\|_{L_2}}{\|\sigma_{true}\|_{L_2}}=0.2552$. The optimal control framework is implemented without regularization $(\beta=0)$ in all but the last subcase (Figure \ref{fig_2D1T_04}(c)). 

\begin{figure}[h!]
\centering
\includegraphics[width=0.9\textwidth]{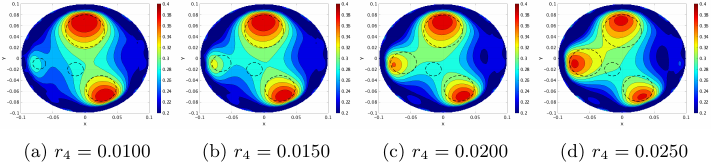}
\caption{Contour plot of obtained conductivity $\sigma_{250}$ for different values of radius $r_4$ in stage 2. Dashed lines represent  size and location of each tumor of true conductivity.}
\label{fig_2D1T_05}
\end{figure}

\vskip.1in
\noindent\textbf{Sensitivity with respect to size.}
For this analysis, we considered different values of $c_4$ and $r_4$. We steadily increased the radius $r_4$ of the smallest tumor while its center $c_3$ was recalculated in order to preserve the distance to the boundary. Figure \ref{fig_2D1T_05} shows the reconstructed conductivity $\sigma_{250}$ for each case of radius $r_4$ and corresponding center $c_4$. Identification of the corresponding tumor cell is improved as the radius $r_4$ increases. Table \ref{tab3} shows the cost value and relative error of voltage and conductivity at the last iteration of stage 2 for each case of radius $r_4$.

\begin{table}[h]
\caption{Metrics for the 2D - Case: 4 tumor - Size.\label{tab3}}
\begin{tabular}{cccc}
\toprule
\multirow{2}{5em}{Radius ($r_4$)} & \multirow{2}{5em}{Cost Value ($\mathcal{K}_{end}$)} & \multicolumn{2}{c}{Relative Error}  \\
& & Voltage ($\frac{|U^{end}-U^\ast|}{|U^\ast|}$) & Conductivity ($\frac{\|\sigma_{end}-\sigma_{true}\|_{L_2}}{\|\sigma_{true}\|_{L_2}}$) \\
\midrule
0.0100 & 7.2617\mbox{e-07} & 0.0602 & 0.2624 \\
0.0150 & 6.5505\mbox{e-07} & 0.0581 & 0.2604 \\
0.0200 & 5.1999\mbox{e-07} & 0.0555 & 0.2516 \\
0.0250 & 6.4377\mbox{e-07} & 0.0504 & 0.2439 \\
\bottomrule
\end{tabular}
\end{table}

\vskip.1in
\noindent\textbf{Regularization effect.}
Here, we have considered the effect of regularization for the convergence of the GPM algorithm. We set initial conditions $\sigma_{ini}$ and $U^{ini}$ as the reconstructed conductivity $\sigma_{end}$ and voltage vector $U^{end}$ obtained from stage 3. Next, we set the regularization parameter $\beta=0.1$ and set $U^\ast$ as the measured voltage obtained from stage 1. The result after 1500 iterations of stage 2 is displayed in Figure \ref{fig_2D1T_04}(c) above. Cost value after 1500 iterations is $\mathcal{K}_{1500}=4.0569\mbox{e-08}$ and relative errors are $\frac{|U_{1500}-U^\ast|}{|U^\ast|}=2.3743\mbox{e-04}$ and $\frac{\|\sigma_{1500}-\sigma_{true}\|_{L_2}}{\|\sigma_{true}\|_{L_2}}=0.1323$, respectively.

\subsection{Results in 3D}  
\begin{figure}[h!]
\centering
\includegraphics[width=0.5\textwidth]{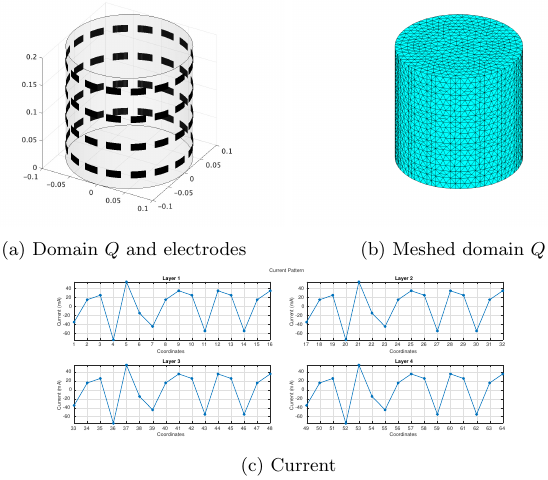}
\caption{(Top Left) Cylindrical domain $Q$ with size and position of electrodes. (Top Right) Meshed domain $Q$. (Bottom) Current pattern used for all cases.}
\label{fig_3D_01}
\end{figure}

In this case, we set $Q$ to be the cylinder of radius $r=0.1\mbox{m}$ and height $h=0.2\mbox{m}$, namely
\begin{equation}
    Q = \{ (x,y,z)\in\R^3:~x^2+y^2 < r^2,\quad 0 < z < h \}.
\end{equation}
A set of $m=64$ electrodes with dimension $0.024~\mbox{rad}$ width and $0.012\mbox{m}$ height were arranged in 4 layers placed in the lateral boundary of $Q$, see Figure \ref{fig_3D_01}(a). These layers are numbered 1--4 from the bottom to the top of the cylinder. Figure \ref{fig_3D_01}(b) shows the linear mesh domain $Q$, consisting of 9392 nodes and 49058 elements (tetrahedrons) with max (edge) size 0.01. A uniform contact impedance vector $Z=(Z_l)_{l=1,\dots,m}$ with $Z_l=0.1~\mbox{Ohm}$ was set. Background conductivity is set to $0.2~(\mbox{Ohm}\cdot m)^{-1}$ corresponding to healthy tissue, and we assume that the conductivity of cancerous tissue is twice as high. Current pattern vector $I$ is set by replicating the pattern used for the 2D cases to each layer of electrodes; see Figure \ref{fig_3D_01}(c).

For the optimization process using the GPM algorithm, we set initial conductivity map $\sigma_{ini}=0.3~(\mbox{Ohm}\cdot\mbox{m})^{-1}$ and initial voltage vector $U^{ini}=(U^{ini}_l)_{l=1}^m$ to $U^{ini}_l=1~\mbox{volt}$ if $l$ is even and $U^{ini}_l=-1~\mbox{volt}$ if $l$ is odd. Additionally, we set termination conditions to a maximum number of iterations $N_{max}$ or the relative error tolerance $\epsilon=10^{-6}$. Finally, for all the cases listed below, we define $Q_\varepsilon = \{(x,y,z)\in~Q:~x^2+y^2<(0.1-\varepsilon)^2\}$ where $\varepsilon>0$ is given.

\subsubsection{Case: 1 Tumor}
\begin{figure}[h!]
\centering
\includegraphics[width=0.7\textwidth]{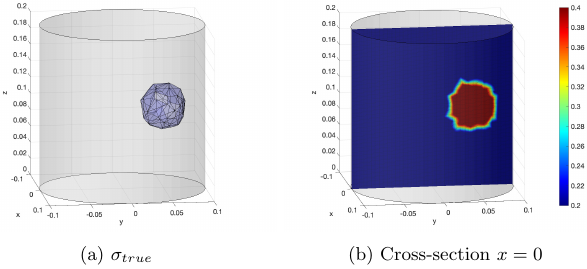}
\caption{(Left) 3D representation of true conductivity tumor showing size and location. (Right) Contour plot of cross-section $x=0$ of true conductivity $\sigma_{true}$.}
\label{fig_3D1T_01}
\end{figure}

Let us consider $\sigma_{true}:Q\to\R$ determining the spherical tumor with center $c=(0,0.05,0.1)$ and radius $r=0.03$. Figure \ref{fig_3D1T_01}(a) shows 3D representation of the conductivity $\sigma_{true}$, its position and size within $Q$, and Figure \ref{fig_3D1T_01}(b) shows the cross-section $x=0$ of $\sigma_{true}$.

\begin{figure}[h!]
\centering
\includegraphics[width=0.7\textwidth]{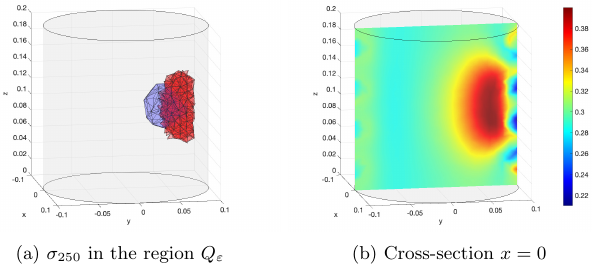}
\caption{(Left) 3D representation of obtained conductivity $\sigma_{250}$ in the region $Q_\varepsilon$. (Right) Contour plot of cross-section $x=0$ of $\sigma_{350}$ for stage 2.}
\label{fig_3D1T_02}
\end{figure}

The optimal control framework is implemented without regularization $(\beta=0)$.
Figure \ref{fig_3D1T_02}(a) shows 3D reconstruction of $\sigma_{250}$ for Stage 2 in the region $\{(x,y,z)\in~Q_\varepsilon : \sigma_{end}(x,y,z)>0.37\}$, for $\varepsilon=10^{-2}$. Figure \ref{fig_3D1T_02}(b) shows the vertical cross-section $x=0$ of the corresponding $\sigma_{250}$.

\begin{figure}[h!]
\centering
\includegraphics[width=0.9\textwidth]{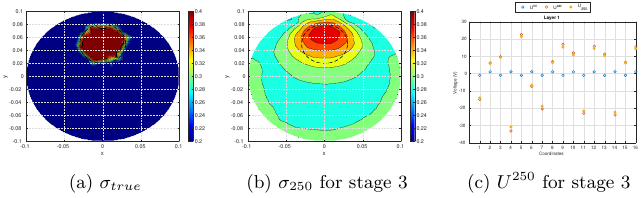}
\caption{Contour plot of the cross-section $z=0.1$ of true conductivity $\sigma_{true}$ (left) and obtained conductivity $\sigma_{250}$ (center) after stage 2. (Right) Coordinates $U_l$, $l=1,\dots,16$, of obtained voltage after stage 3.}
\label{fig_3D1T_03}
\end{figure}

Figure \ref{fig_3D1T_03}(a)-(b) shows the contour plot of the cross-section $z=0.1$ (vertical center of $Q$) of the true conductivity map $\sigma_{true}$ and calculated conductivity $\sigma_{end}$ at the last iteration for stage 2, respectively. Dashed lines correspond to the position of the true conductivity $\sigma_{true}$. 
Figure \ref{fig_3D1T_03}(c) shows the coordinates $l=1,\dots,16$ (first layer of electrodes) of the obtained voltage vector $U^{end}$ for stage 2 in the 1 Tumor case, voltage values at the remaining electrodes match those for the first layer ($l=1,\dots,16$). Cost value at the final iteration of stage 3 is $\mathcal{K}_{end}=7.7029\mbox{e-08}$ and relative errors of voltage and conductivity are $\frac{|U^{end}-U^\ast|}{|U^\ast|}=0.0696$ and $\frac{\|\sigma_{end}-\sigma_{true}\|_{L_2}}{\|\sigma_{true}\|_{L_2}}=0.4875$

\begin{figure}[h!]
\centering
\includegraphics[width=0.9\textwidth]{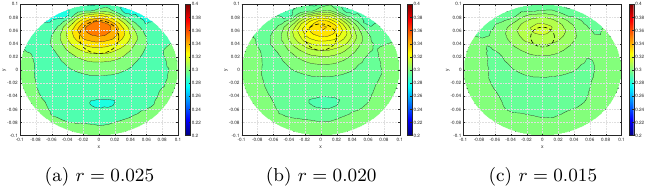}
\caption{Contour plot of the cross-section $z=0.1$ of obtained conductivity $\sigma_{250}$ for different values of radius $r$ for $\sigma_{true}$.}
\label{fig_3D1T_04}
\end{figure}

\vskip.1in
\noindent\textbf{Sensitivity with respect to size.} For this analysis, we have considered $\sigma_{true}$ for fixed center $c=(0,0.05,0.1)$ and three different values of radius, namely $r=0.03,0.02,0.01$. Figure \ref{fig_3D1T_04} shows the cross-section $z=0.1$ of the reconstructed conductivity $\sigma_{N}$ and the last iteration $N=250$ for the cases $r=0.025,0.020,0.015$. Dashed lines show the location and size of $\sigma_{true}$. Table \ref{tab4} shows the values of cost functional end relative error with respect to voltage and conductivity at the last iteration of stage 2 for each case of radius $r$.

\begin{table}[h]
\caption{Metrics for the 3D - Case: 1 tumor - Size.\label{tab4}}
\begin{tabular}{cccc}
\toprule
\multirow{2}{5em}{Radius ($r$)} & \multirow{2}{5em}{Cost Value ($\mathcal{K}_{end}$)} & \multicolumn{2}{c}{Relative Error}  \\
& & Voltage ($\frac{|U^{end}-U^\ast|}{|U^\ast|}$) & Conductivity ($\frac{\|\sigma_{end}-\sigma_{true}\|_{L_2}}{\|\sigma_{true}\|_{L_2}}$) \\\midrule
0.030 & 7.7029\mbox{e-08} & 0.0697 & 0.4876 \\
0.025 & 2.7787\mbox{e-08} & 0.0698 & 0.4884 \\
0.020 & 1.2547\mbox{e-08} & 0.0699 & 0.4974 \\
0.015 & 3.6822\mbox{e-09} & 0.0700 & 0.5008 \\
0.010 & 2.1157\mbox{e-09} & 0.0700 & 0.5023 \\
\bottomrule
\end{tabular}
\end{table}

\begin{figure}[h!]
\centering
\includegraphics[width=0.9\textwidth]{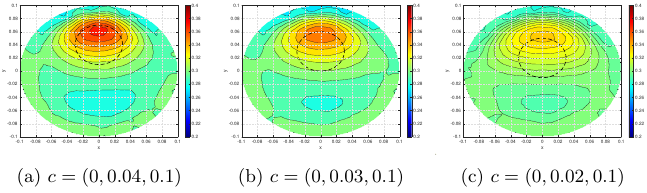}
\caption{Contour plot of the cross-section $z=0.1$ of obtained conductivity $\sigma_{250}$ for different values of center $c$ for $\sigma_{true}$.}
\label{fig_3D1T_05}
\end{figure}

\vskip.1in
\noindent\textbf{Sensitivity with respect to location.}
Let us considered $\sigma_{true}$ with fixed radius $r=0.03$ and different center positions, namely $c = (0,0.05,0.1),(0,0.03,0.1),(0,0.01,0.1)$. Figure \ref{fig_3D1T_05} shows the cross-section contour plot at $z=0.1$ of obtained conductivity for the cases $c=(0,y,0.1)$, with $y=0.04,0.03,0.02$. Table \ref{tab5} shows the corresponding cost functional values and relative error with respect to voltage and conductivity at the final iteration of stage 2 for each case of center $c$.

\begin{table}[h]
\caption{Metrics for the 3D - Case: 1 tumor - Location.\label{tab5}}
\begin{tabular}{cccc}
\toprule
\multirow{2}{7em}{Center ($c=(0,y,0.1)$)} & \multirow{2}{5em}{Cost Value ($\mathcal{K}_{end}$)} & \multicolumn{2}{c}{Relative Error}  \\
& & Voltage ($\frac{|U^{end}-U^\ast|}{|U^\ast|}$) & Conductivity ($\frac{\|\sigma_{end}-\sigma_{true}\|_{L_2}}{\|\sigma_{true}\|_{L_2}}$) \\
\midrule
y=0.05 & 7.7029\mbox{e-08} & 0.0697 & 0.4876 \\
y=0.04 & 6.5934\mbox{e-07} & 0.0697 & 0.4895 \\
y=0.03 & 2.9617\mbox{e-08} & 0.0698 & 0.4850 \\
y=0.02 & 2.5079\mbox{e-08} & 0.0698 & 0.4862 \\
y=0.01 & 2.0786\mbox{e-08} & 0.0698 & 0.4870 \\
y=0.00 & 1.9150\mbox{e-08} & 0.0698 & 0.4869 \\
\bottomrule
\end{tabular}
\end{table}

\subsubsection{Case: 2 Tumors}
\begin{figure}[h!]
\centering
\includegraphics[width=0.9\textwidth]{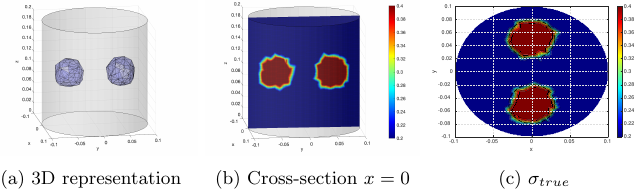}
\caption{(Left) 3D representation of the tumor of the true conductivity $\sigma_{true}$. (Center) contour plot of the cross-section $x=0$ of $\sigma_{true}$. (Right) Contour plot of the cross-section $z=0.1$ of $\sigma_{true}$.}
\label{fig_3D2T_01}
\end{figure}

In this section, we consider $\sigma_{true}:Q\to\mathbb{R}$ determining two spherical tumors with center $c_1=(0,0.05,0.1)$, $c_2=(0,-0.05,0.1)$ and radius $r_1=r_2=0.03$. 
Figures \ref{fig_3D2T_01} show the 3D representation of $\sigma_{true}$, the vertical cross-section $x=0$, and the horizontal cross-section $z=0.1$ (center of the cylinder) of conductivity $\sigma_{true}$, respectively. 

\begin{figure}[h!]
\centering
\includegraphics[width=0.9\textwidth]{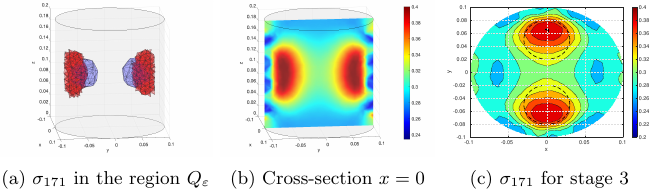}
\caption{(Left) 3D representation of obtained conductivity $\sigma_{171}$ for stage 2 in the region $Q_\varepsilon$. (Center) contour plot of the cross-section $x=0$ of $\sigma_{171}$. (Right) Contour plot of the cross-section $z=0.1$ of $\sigma_{171}$.}
\label{fig_3D2T_02}
\end{figure}

The optimal control framework is implemented without regularization $(\beta=0)$.
Figure \ref{fig_3D2T_02}(a) shows the reconstruction of calculated conductivity $\sigma_{171}$ for stage 3 in the region $\{\mathbf{x} \in Q_\varepsilon : \sigma_{171}(\mathbf{x}) > 0.37\}$, for $\varepsilon=10^{-2}$. Figures \ref{fig_3D2T_02}(b)-(c) show the cross-sections $x=0$ and cross-section $z=0.1$ of $\sigma_{171}$. Dashed lines show the position and size of true conductivity $\sigma_{true}$. 

Finally, the cost values at the last iteration of stage 2 is $\mathcal{K}_{end}=7.7029\mbox{e-08}$ and corresponding relative errors for voltage and conductivity are $\frac{|U^{end}-U^\ast|}{|U^\ast|}=0.0683$ and $\frac{\|\sigma_{end}-\sigma_{true}\|_{L_2}}{\|\sigma_{true}\|_{L_2}}=0.4527$.

\begin{figure}[h!]
\centering
\includegraphics[width=0.9\textwidth]{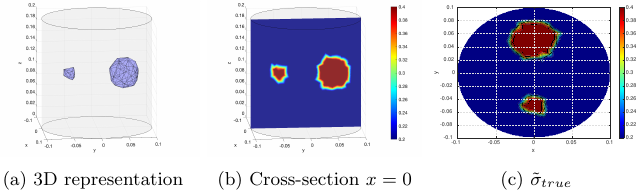}
\caption{(Left) 3D representation of the tumors of the true conductivity $\sigma_{true}$. (Center) contour plot of cross-section $x=0$ of $\sigma_{true}$. (Right) Contour plot of the cross-section $z=0.1$ of $\sigma_{true}$.}
\label{fig_3D2T_03}
\end{figure}

\vskip.1in
\noindent\textbf{Sensitivity with respect to size.}
In this case, we consider the case in which the size of one of the tumors determined by initial true conductivity $\sigma_{true}$ is reduced. Indeed, we set $r_2=0.015$, positions $c_1$, $c_2$ and size $r_1$ are kept the same. The obtained new conductivity will be denoted $\tilde{\sigma}_{true}$. Figures \ref{fig_3D2T_03} show the 3D representation of $\tilde{\sigma}_{true}$, and corresponding cross-sections $x=0$ and $z=0.1$, respectively.

\begin{figure}[h!]
\centering
\includegraphics[width=0.9\textwidth]{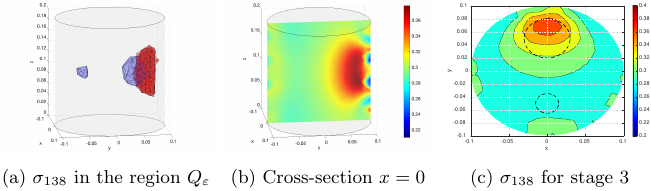}
\caption{(Left) 3D representation of the obtained conductivity $\sigma_{138}$ in the region $Q_\varepsilon$ after stage 2. (Center) contour plot of cross-section $x=0$ of $\sigma_{138}$. (Right) Contour plot of the cross-section $z=0.1$ of $\sigma_{138}$.}
\label{fig_3D2T_04}
\end{figure}

The optimal control framework is implemented without regularization $(\beta=0)$. 
Figure \ref{fig_3D2T_04}(a) shows the reconstruction of calculated conductivity $\sigma_{138}$ at the last iteration for stage 2 in the region $\{\mathbf{x}\in Q_\varepsilon : \sigma_{138}(\mathbf{x})>0.35\}$, for $\varepsilon=10^{-2}$. Figures \ref{fig_3D2T_04}(b)-(c) show cross-sections $x=0$ and $z=0$ of $\sigma_{138}$, respectively. Dashed lines show the position and size of true conductivity $\tilde{\sigma}_{true}$.

\begin{figure}[h!]
\centering
\includegraphics[width=0.6\textwidth]{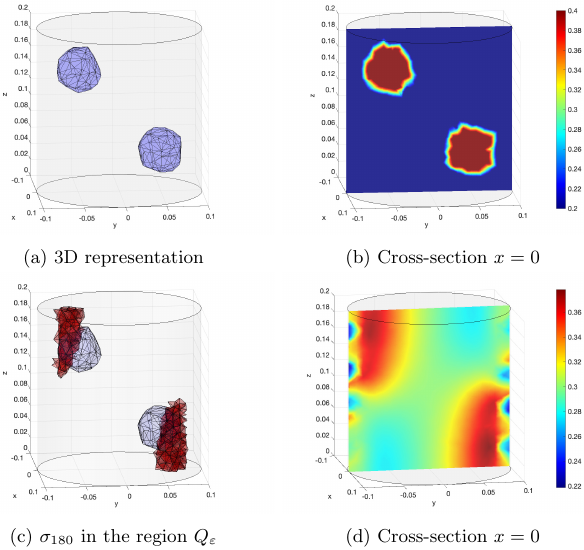}
\caption{(Top Left) 3D representation of the true conductivity $\sigma_{true}$. (Top Right) Contour plot of the cross-section $x=0$ of $\sigma_{true}$. (Bottom Left) Obtained conductivity $\sigma_{180}$ in the region $Q_\varepsilon$ after stage 2. (Bottom Right) contour plot of the cross section $x=0$ of $\sigma_{180}$.}
\label{fig_3D2T_05}
\end{figure}

\vskip.1in
\noindent\textbf{Sensitivity with respect to location.}
In this case, we consider the case in which the center of the tumors determined by the initial true conductivity $\sigma_{true}$ are modified. Indeed, we set $c_1=(0,0.05,0.05)$, $c_2=(0,-0.05,0.15)$, the radius $r_1=r_2=0.03$ are preserved. The new conductivity will be denoted $\bar{\sigma}_{true}$. Figure \ref{fig_3D2T_05}(a)-(b) show the 3D representation of $\bar{\sigma}_{true}$ and cross-section $x=0$, respectively.
The optimal control framework is implemented without regularization $(\beta=0)$. 
Figure \ref{fig_3D2T_05}(c) shows the reconstructed $\sigma_{180}$ for stage 2 at the last iteration in the region $\{\mathbf{x}\in Q_\varepsilon : \sigma_{180}(\mathbf{x})>0.35\}$, for $\varepsilon=10^{-2}$, while Figure \ref{fig_3D2T_05}(d) shows the cross-section $x=0$ of corresponding conductivity $\sigma_{180}$.

\begin{figure}[h!]
\centering
\includegraphics[width=0.6\textwidth]{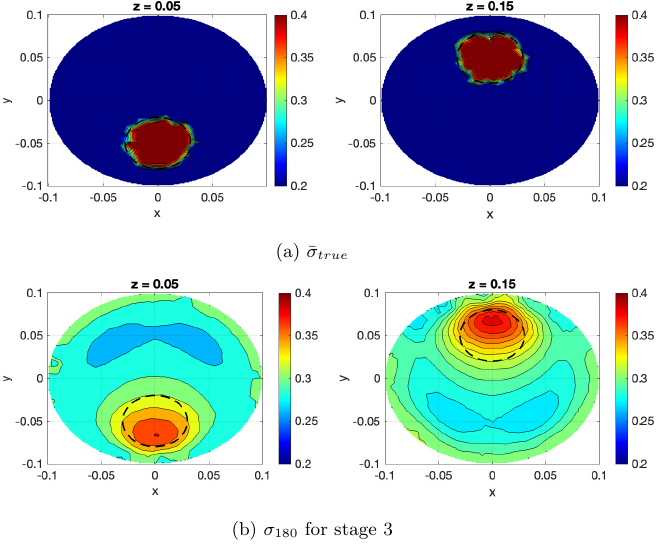}
\caption{(Top) Contour plot of cross-sections $z=0.05$ and $z=0.15$ of the true conductivity $\sigma_{true}$.(Bottom) Contour plot of cross-sections $z=0.05$ and $z=0.15$ of obtained conductivity $\sigma_{180}$ after stage 2.}
\label{fig_3D2T_06}
\end{figure}

Finally, Figures \ref{fig_3D2T_06} show cross-sections $z=0.05$ and $z=0.15$, across the center of each tumor, of true conductivity $\bar{\sigma}_{true}$ and obtained conductivity $\sigma_{180}$ for stage 2 at the last iteration, respectively.

\subsubsection{Case: 4 Tumors}
\begin{figure}[h!]
\centering
\includegraphics[width=0.9\textwidth]{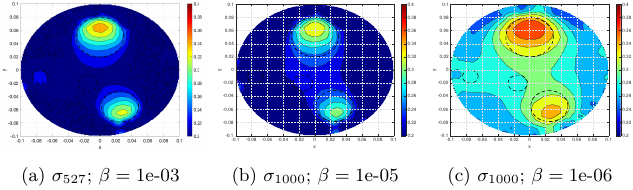}
\caption{Contour plot of the cross-section $z=0.1$ of the obtained conductivity $\sigma_{end}$ after stage 2 for different values of regularization parameter $\beta$.}
\label{fig_3D4T_04}
\end{figure}

In this section, we consider the case of $\sigma_{true}$ determining four spherical tumors of center $c_1=(0,0.05,0.1)$, $c_2=(-0.075,-0.01,0.1)$, $c_3=(-0.015,-0.02,0.1)$, $c_4=(0.025,-0.055,0.1)$, and respective radius $r_4=0.03$, $r_2=0.0099$, $r_3=0.15$ and $r_4=0.02$.
The optimal control framework is implemented without regularization $(\beta=0)$ in all but the last subcase (Figure \ref{fig_3D4T_04}). 

\begin{figure}[h!]
\centering
\includegraphics[width=0.7\textwidth]{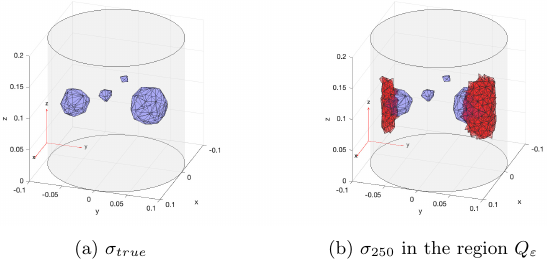}
\caption{(Left) 3D representation of the tumors of the true conductivity $\sigma_{true}$. (Right) Obtained conductivity $\sigma_{250}$ in the region $Q_\varepsilon$ after stage 2.}
\label{fig_3D4T_01}
\end{figure}

Figure \ref{fig_3D4T_01} shows 3D representation of $\sigma_{true}$ and reconstructed conductivity $\sigma_{250}$ within the region $\{\mathbf{x}\in~Q_\varepsilon:~\sigma_{end}(\mathbf{x})>0.35\}$, for $\varepsilon=10^{-2}$, respectively. 

\begin{figure}[h!]
\centering
\includegraphics[width=0.7\textwidth]{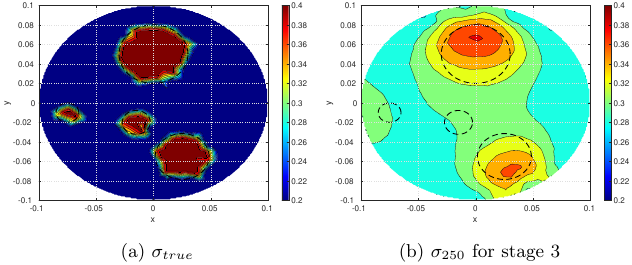}
\caption{Contour plot of the cross-section $z=0.1$ of the true conductivity $\sigma_{true}$ (Left) and obtained conductivity $\sigma_{250}$ (Right) after stage 2.}
\label{fig_3D4T_02}
\end{figure}

Figure \ref{fig_3D4T_02} shows the cross-section $z=0.1$ (center of the cylinder) of $\sigma_{true}$ and obtained conductivity $\sigma_{end}$ at the last iteration for stage 2. Cost value at the last iteration of stage 2 is $\mathcal{K}_{end}=7.7029\mbox{e-08}$ and relative errors for voltage and conductivity with respect to measured voltage $U^\ast$ and true conductivity $\sigma_{true}$ are $\frac{|U^{end}-U^\ast|}{|U^\ast|}=0.0697$ and $\frac{\|\sigma_{end}-\sigma_{true}\|_{L_2}}{\|\sigma_{true}\|_{L_2}}=0.4876$.

\begin{figure}[h!]
\centering
\includegraphics[width=0.9\textwidth]{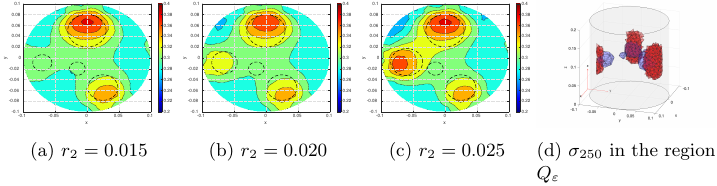}
\caption{Contour plot of the cross-section $z=0.1$ of the obtained conductivity $\sigma_{250}$ after stage 2 for different values of radius $r_2$.}
\label{fig_3D4T_03}
\end{figure}

\vskip.1in
\noindent\textbf{Sensitivity with respect to size.}
Here we considered different values for the tumor cell with center $c_2$ and radius $r_2$ of $\sigma_{true}$. Indeed, the radius $r_2$ is increased while the center $c_2$ is recalculated in order to preserve the distance to the lateral boundary of $Q$. Figures \ref{fig_3D4T_03}(a)-(c) show the cross-section $z=0.1$ (center of the cylinder) of $\sigma_{end}$ for all the cases of radius $r_2$. Figure \ref{fig_3D4T_03}(d) shows 3D reconstruction of the region $\{\mathbf{x}\in Q_\varepsilon:\sigma_{250}(\mathbf{x})>0.35\}$, for $\varepsilon=10^{-2}$. Table \ref{tab6} shows the cost value and relative errors of voltage and conductivity at the last iteration of stage 2 and for each case radius $r_2$.

\begin{table}[h]
\caption{Metrics for the 3D - Case: 4 tumors - Size.\label{tab6}}
\begin{tabular}{cccc}
\toprule
\multirow{2}{5em}{Radius ($r_2$)} & \multirow{2}{5em}{Cost Value ($\mathcal{K}_{end}$)} & \multicolumn{2}{c}{Relative Error}  \\
& & Voltage ($\frac{|U^{end}-U^\ast|}{|U^\ast|}$) & Conductivity ($\frac{\|\sigma_{end}-\sigma_{true}\|_{L_2}}{\|\sigma_{true}\|_{L_2}}$) \\\midrule
0.015 & 9.3826\mbox{e-08} & 0.0693 & 0.4703 \\
0.020 & 1.0858\mbox{e-07} & 0.0691 & 0.4632 \\
0.025 & 8.8688\mbox{e-08} & 0.0692 & 0.4600 \\
\bottomrule
\end{tabular}
\end{table}

\vskip.1in
\noindent\textbf{Regularization effect.}
We have considered here the effect of regularization in the reconstruction/optimization process. Initial conditions $\sigma_{ini}$ and $U^{ini}$ were set to those obtained after 250 iterations in stage 2 without regularization (see Figure \ref{fig_3D4T_02}(b)). Figure \ref{fig_3D4T_04} shows the cross section contour plot of resulting conductivity (with values in $(\mbox{Ohm}\cdot\mbox{m})^{-1}$) for different values of regularization parameter $\beta$. Table \ref{tab7} shows cost value and relative errors of voltage and conductivity at the last iteration. The minimum value of the cost functional, and relative error of control parameters are minimized at the value $10^{-3}$ of the regularization parameter.  

\begin{table}[]
\caption{Metrics for the 3D - Case: 4 tumors - Regularization.\label{tab7}}
\begin{tabular}{cccc}
\toprule
\multirow{2}{5em}{Parameter ($\beta$)} & \multirow{2}{5em}{Cost Value ($\mathcal{K}_{end}$)} & \multicolumn{2}{c}{Relative Error}  \\
& & Voltage ($\frac{|U^{end}-U^\ast|}{|U^\ast|}$) & Conductivity ($\frac{\|\sigma_{end}-\sigma_{true}\|_{L_2}}{\|\sigma_{true}\|_{L_2}}$) \\\midrule
1\mbox{e-03} & 4.2204\mbox{e-08} & 9.3427\mbox{e-06} & 0.0910 \\
1\mbox{e-05} & 1.7765\mbox{e-07} & 9.0933\mbox{e-04} & 0.1027 \\
1\mbox{e-06} & 6.9503\mbox{e-06} & 0.0191 & 0.3554 \\
\bottomrule
\end{tabular}
\end{table}

\section*{Discussion}\label{sect:Dis}

We consider an inverse EIT problem on the identification of the conductivity map in the complete electrode model based on the $m$ current-to-voltage measurements on the boundary electrodes. Particular motivation arises from the medical application for the identification of the cancerous tumor at early stages of development. The idea of the method is based on the fact that the electrical conductivity of the cancerous tumor is significantly higher than the conductivity of the healthy tissue. A variational formulation as a PDE-constrained optimal control problem is introduced based on the novel idea of increasing the size of the input data by adding "voltage-to-current" measurements through various permutations of the single "current-to-voltage" measurement. The idea of permutation preserves the size of the unknown parameters at the expense of increasing the number of PDE constraints. We apply a Gradient Projection Method based on the Fr\'echet differentiability in Besov-Hilbert spaces. 
\begin{itemize}
\item[-] Numerical simulations demonstrate that for both 2D and 3D model examples, the resolution of target tumor regions is significantly improved by increasing the number of input data from $m$ to $m^2$.
\item[-] Resolution of target tumor regions is demonstrated to be sensitive to the size of the tumor and its distance from the boundary electrodes. Smaller tumor size and greater distance from the electrodes negatively impact the resolution of tumors produced by the method. 
\item[-] Based on the effective computational performance a new 2-step procedure is suggested for the medical application for the identification of the cancerous tumor at an early stage of its development. 
\end{itemize}

\section*{Conclusions}\label{sect:Dis}

This paper suggests a new method for the identification of the cancerous tumor at an early stage of development. Relying on the experimental fact that the electrical conductivity of the cancerous tumor is significantly higher than the conductivity of the healthy tissue, we consider
an inverse EIT problem on the identification of the conductivity map in the complete electrode model based on the $m$ current-to-voltage measurements on the boundary electrodes.  
A variational formulation as a PDE-constrained optimal control problem is introduced. To address the ill-posedness of the inverse problem due to insufficient measurements, we implement a novel idea of increasing the size of the input data by adding "voltage-to-current" measurements through various permutations of the single "current-to-voltage" measurement. The idea of permutation preserves the size of the unknown parameters on the expense of increase of the number of PDE constraints. We apply a gradient projection method based on the Fr\'echet differentiability in Besov-Hilbert spaces. 
Numerical simulations of 2D and 3D model examples demonstrate the sharp increase of the resolution of the cancerous tumor by increasing the number of measurements from $m$ to $m^2$. Based on the effective computational performance a new 2-step procedure is suggested for the identification of the cancerous tumor at an early stage of its development in the clinical setting. 

\section*{Ethics approval and comsent to participate}
Not applicable
\subsection*{Consent for publication}
Not applicable.
\subsection*{Availability of data and materials}
The datasets generated and/or analysed during the current study are not publicly available due to potential conflicts with other research groups but are available from the corresponding author on reasonable request.
\subsection*{Competing interests}
All other authors declare they have no competing interests.
\subsection*{Author contributions}
UGA initiated the idea of the new method, pursued all the mathematical proofs and rigorous analysis, and supervised the project; JHR pursued coding and software development for the numerical implementation of the new method; UGA and JHR pursued validation and computational analysis of numerical simulations; UGA drafted the manuscript; Both UGA and JHR read and agreed with the final version of the manuscript; 
\subsection*{Funding}
No funding was received.
\subsection*{Ethics approval and consent to participate}
Not applicable.
\subsection*{Acknowledgements} 
Authors are grateful for the help and support provided by the Scientific Computing and Data Analysis section of Core Facilities at OIST.

\newpage
\bibliography{ref.bib}

\end{document}